\documentstyle[amssymb]{amsart} 

\newcommand{\ZFCs}{\text{\normalshape\sf ZFC$^\star$}}
\newcommand{\ZFCa}{\text{\normalshape\sf ZFC}}

\newcommand{\reals}{\Bbb R}
\newcommand{\rationals}{\Bbb Q}
\newcommand{\rest}{{\restriction}}
\newcommand{\add}{\text{\normalshape\sf{add}}}
\newcommand{\cov}{\text{\normalshape\sf{cov}}}
\newcommand{\unif}{\text{\normalshape\sf{non}}}
\newcommand{\cof}{\text{\normalshape\sf{cof}}}
\newcommand{\cf}{{\text{\normalshape\sf{cf}}}}

\newcommand{\suc}{{\text{\normalshape\sf {succ}}}}

\newcommand{\QED}{\hspace{0.1in} \Box \vspace{0.1in}}

\newcommand{\N}{{\cal N}}
\newcommand{\M}{{\cal M}}
\newcommand{\V}{{\bold V}}
\newcommand{\<}{\langle}
\renewcommand{\>}{\rangle}

\newcommand{\thinks}{\models}

\newtheorem{theorem}{Theorem}[section]
\newtheorem{lemma}[theorem]{Lemma}

\newtheorem{definition}[theorem]{Definition}

\newcommand{\lesdot}{\mathrel{\mathord{<}\!\!\raise 
0.8 pt\hbox{$\scriptstyle\circ$}}}
\newcommand{\Proof}{{\sc Proof} \hspace{0.2in}}

\newcommand{\lft}[2]{\mathopen\ifcase#1{}\oo\or
                        \big#2\or\Big#2\else\oo\fi} 
\newcommand{\rgt}[2]{\mathclose\ifcase#1{}\oo\or
                        \big#2\or\Big#2\else\oo\fi} 

\newcommand{\R}{\text{\normalshape\sf R}}
\begin{document}

\title{Borel images of sets of reals}
\author{Tomek Bartoszy\'{n}ski}
\address{Department of Mathematics\\
Boise State University\\
Boise, Idaho 83725, USA}
\email{{\tt tomek@@math.idbsu.edu}}
\author{Haim Judah}
\address{Abraham Fraenkel Group for Mathematical Logic\\
Department of Mathematics\\
Bar Ilan University\\
52900 Ramat Gan, Israel 
}
\email{{\tt judah@@bimacs.cs.biu.ac.il}}
\thanks{Research partially supported by 
the Israeli Academy of Science, Basic Research Foundation and FONDECYT
\#1940695, Chile}
\maketitle
\begin{abstract}
The main goal of this paper is to generalize several results
concerning cardinal invariants to the statements about the 
associated families of sets. We also discuss the relationship between
the additive properties of sets and their Borel images. 
Finally, we present  estimates for the size of the smallest set
which is not strongly meager.
\end{abstract} 
\section{Introduction}
The purpose of this paper is to study additive properties of sets of
reals. By reals we mean the Cantor set $2^\omega$, real line $\reals$
or interval $[0,1]$.
We will be working in the space $2^\omega$ with addition modulo 2, but most of
the results would translate to the interval $[0,1]$ and $\reals$.

We will use the standard notation. Let $\rationals$ be the canonical
countable 
dense set in the spaces mentioned above. It is the set of rationals in
case of $\reals $ or $[0,1]$ and a collection of 0-1 sequences that
are eventually equal to zero in case of $2^\omega$. 

We will often be using trees on $2^{<\omega}$. A subset $T \subseteq
2^{<\omega}$ is a tree if for every $t \in T$ and $n <|t|$, $t
\rest n \in T$. We will also require that $T$ does not have terminal
nodes.
Let $T \rest n$ denote the $n$-th level of $T$ and for $n<m$ and $s \in
T \rest n$ let 
$$\suc_{T,m}(s) = \{t \in T \rest m : s \subseteq t\}.$$

For a tree $T$ and $s \in T$ let $T_s = \{t \in T : s \subseteq t
\text{ or } t \subseteq s\}$ be the subtree determined by $s$.

For a set $A \subseteq 2^{<\omega}$ let 
$$[A] = \{x \in 2^\omega: \forall n \ x \rest n \in A \text{ or } x
\text{ contains a terminal node from $A$} \}.$$

Note that in case of a tree, $[T]$ coincides with the set of branches
of $T$. On the other hand, for $s \in 2^{<\omega}$, $[s]$ is the
basic open set in $2^\omega$ determined by $s$.

If $T \subseteq 2^{<\omega}$ is a tree as above 
and $n \in \omega$ we define the tree
$T^{\<n\>}$ as
$$s \in T^{\<n\>} \iff  |s|<n \text{ or } \exists t \in T \
\lft1( |s|=|t| \ \& \ s \rest [n, |s|) = t \rest [n, |s|) \rgt1).$$
Note that if we identify the set of rationals $\rationals$ with
elements of $2^\omega$ which are eventually equal to $0$, then
$[T]+\rationals = \bigcup_{n \in \omega} [T^{\<n\>}]$.

$\ZFCs$ always denotes a finite, sufficiently large fragment of
$\ZFCa$. 

Quantifiers $\forall^\infty$ and $\exists^\infty$ denote ``for
all except finitely many'' and ``for infinitely many'', respectively.

For a set $H \subseteq 2^\omega
\times 2^\omega$ and $x,y \in 2^\omega$ let $(H)_x = \{y : \<x,y\> \in
H\}$ and $(H)^y = \{ x : \<x , y\> \in H\}$. By $\M$ and $\N$ we
denote the $\sigma$-ideals of meager and Lebesgue measure zero sets
respectively. 
\begin{definition}
Suppose that ${\cal J}$ is an ideal of subsets of the real line.
A Borel set $H \subseteq 2^\omega \times 2^\omega$ is called  a
${\cal J}$-set  if $(H)_x \in {\cal J} $ for all $ x \in 2^\omega$. 
 
We say that $X \subseteq 2^\omega $ is an $R^{{\cal J}}$ set if 
for every ${\cal J}$-set $H$,
$$\bigcup_{x \in X} (H)_x \neq 2^\omega .$$
A set $X \subseteq 2^\omega $ is an $SR^{{\cal J}}$ set if 
for every ${\cal J}$-set $H$,
$$\bigcup_{x \in X} (H)_x \in {\cal J} .$$

Consider the space $2^\omega$ with group operation $+$ defined as
addition modulo $2$.

A set $X$ is ${\cal J}$-additive if for every set $F \in {\cal J}$,
$X+F \in {\cal J}$ and is
not ${\cal J}$-covering  if for every $F \in {\cal J}$, $X+F
\neq 2^\omega$. 

Let ${\cal J}^\star$ denote the ideal of ${\cal
  J}$-additive sets. 
\end{definition}
Traditionally the sets which are not $\N$-covering are called strongly
meager while the sets which are not $\M$-covering are called strong measure
zero sets. Let $\cal {SZ}$ denote the ideal of strongly measure zero
sets and $\cal {SM}$ the collection of strongly meager sets.

It is well known that strong measure zero sets form an ideal whereas
it is not clear whether the collection strongly meager sets is an ideal.
 
We have the following easy observation:
\begin{lemma}\label{1.2}
  Suppose that ${\cal J}$ is a proper, translation invariant, ideal
  with Borel basis that
contains all singletons. Then we have the following inclusions
(writing $\rightarrow$ for $\subseteq$):
$$\begin{array}{ccccc}
 \ & \ & \text{not ${\cal J}$-covering} & \ & \  \\
 \ & \nearrow & \ & \nwarrow &   \ \\
 R^{{\cal J}} & \ & \ & \ &  {\cal J}^\star \\
\ & \nwarrow & \ & \nearrow & \ \\
\ & \ & SR^{{\cal J}} & \ & \  \end{array}.$$
In particular,
$$\begin{array}{cccccccccccc}
 \ & \ & \cal {SZ} & \ & \  & \hspace{0.5in}  & \ & \ & \cal {SM} & \ & \ \\
 \ & \nearrow & \ & \nwarrow &   \ & \hspace{0.5in} & \ &  \nearrow & \ &
   \nwarrow & \ \\ 
 R^{\M} & \ & \ & \ &  {\M}^\star & \hspace{0.5in}  & R^{\N} &\ &\ & \ &
 \  {\N}^\star \\ 
\ & \nwarrow & \ & \nearrow & \ & \ & \ & \nwarrow & \ & \nearrow & \\ 
\ & \ & SR^{\M} & \ & \ & \ & \ & \ & SR^{\N} & \ & \end{array} .~\QED$$

\end{lemma}

\begin{definition}
 For any proper ideal ${\cal J}$ of subsets of $X$ we can define the
 following cardinal 
coefficients:

$\add({\cal J}) = \min\left\{ |{\cal A}| : {\cal A} \subseteq {\cal J},
\hbox{ and } \bigcup {\cal A} \not \in {\cal J}\right\}$,

$\cov({\cal J}) = \min\left\{ |{\cal A}| : {\cal A} \subseteq {\cal J},
\hbox{ and } \bigcup {\cal A} = X\right\}$,

$\unif({\cal J}) = \min\left\{ |Y| : Y \subseteq X,
\hbox{ and } Y \not \in {\cal J}\right\}$,

$\cof({\cal J}) = \min\left\{|{\cal A}| : {\cal A} \subseteq {\cal J},
\hbox{ and } \forall B \in {\cal J} \ \exists A \in {\cal A} \ B
\subseteq A \right\}  .$
\end{definition}

Following \cite{vdo:inttop},
let ${\frak b}$ and ${\frak d}$ denote the sizes of the smallest
unbounded and dominating families in $\omega^\omega$ respectively.
 
We have the following lemma:
\begin{lemma}\label{clear}
  If $X$ is an $R^{{\cal J}}$ set then every Borel image of $X$ is an
  $R^{{\cal J}}$ set. 

  If $X$ is an  $SR^{{\cal J}}$ set then every Borel image of $X$
  is an
  $SR^{{\cal J}}$ set. 
\end{lemma}
\Proof
Suppose that $X \subseteq 2^\omega$ and let $f : 2^\omega
\longrightarrow 2^\omega $ be a Borel function.

Let $H \subseteq 2^\omega \times 2^\omega$ be a Borel ${\cal J}$-set.
Define
$$\widetilde{H} = \{\<x,y\> : \<f(x),y\> \in H\} .$$
It is easy to see that $\widetilde{H}$ is a Borel set and that 
$$\bigcup_{x \in X} (H)_{f(x)} = \bigcup_{x \in X} \widetilde{H}_x
.~\QED $$

Note that we have the following easy observation:
\begin{lemma}\label{triv}
  $X \subseteq 2^\omega $ is an $R^{\cal J}$-set iff for every Borel set $H
  \subseteq 2^\omega \times 2^\omega $, such that
  $(H)_x \in {\cal J} $ for 
  all $x \in 2^\omega$, 
$$2^\omega \setminus \bigcup_{x \in X} (H)_x \not \in {\cal J}. $$
\end{lemma}
\Proof Implication $(\leftarrow)$ is obvious. 

On the other hand
suppose that $G$ is a Borel set such that 
$2^\omega \setminus \bigcup_{x \in X} (H)_x \subseteq G
\in {\cal J}$. Let $\widehat{H} = H 
\cup (2^\omega \times G)$. Clearly $\widehat{H}$ witnesses that $X$ is
not an $R^{\cal J}$-set.~$\QED$

Note that
\begin{theorem}[\cite{Kur66Top} p. 434]
  Suppose that $X \subseteq 2^\omega$. 

If $H \subseteq X \times 2^\omega$ is a Borel set then there exists a
Borel set $\widetilde{H} \subseteq 2^\omega \times 2^\omega$ such that
$H = \widetilde{H} \cap (X \times 2^\omega)$.

If $f : X \longrightarrow 2^\omega$ is a Borel function then there
exists a Borel function $\widetilde{f} : 2^\omega \longrightarrow
2^\omega$ such that $\widetilde{f} = f {\rest} X$.~$\QED$
\end{theorem}

Finally we will need the following theorem concerning representation
of $\M$- and $\N$-sets.

\begin{lemma}[Fremlin]\label{fremlin}
Suppose that $H \subseteq 2^\omega \times 2^\omega$ is a Borel set.
\begin{enumerate}
\item Assume $(H)_x$ is meager for all $x$. Then there
exists a sequence of Borel sets $\{G_n :n \in \omega\} \subseteq
2^\omega \times 2^\omega $ such that
\begin{enumerate}
  \item $(G_n)_x$ is a closed nowhere dense set for all $x \in
2^\omega$,
\item $H \subseteq \bigcup_{n \in \omega} G_n$.
\end{enumerate}
\item For every $\varepsilon>0$
there exists a Borel set $B \subseteq 2^\omega \times
2^\omega$ such that 
\begin{enumerate}
\item $H \subseteq B$,
\item $(B)_x$ is open for every $x$,
\item $\mu((B \setminus H)_x) < \varepsilon$ for every $x$.
\end{enumerate}
\end{enumerate}
  \end{lemma}
\Proof
For completeness we present a sketch of the proof here.

Let  $\cal G$  be the family of Borel subsets  $G$  of $2^\omega \times
2^\omega$  such that $(G)_x$ is open for every $x \in 2^\omega$.

(1) 
Let ${\cal J}$ be the $\sigma$-ideal of subsets of the plane generated
by Borel sets $F$ such that $(F)_x$ is closed nowhere dense for all $x$.

Consider the family $\Sigma$ of subsets  $E$  of the plane
such that there are two Borel sets  $G$, $H$  of which  $G \in\cal G$,  
$H \in {\cal J}$ and
$E\triangle G\subseteq H$.   
By induction on Borel class we will show that $\Sigma$ contains all
Borel sets. In particular, it follows that all $\M$-sets are in ${\cal
  J}$, which will finish the proof.
 
Clearly  $\Sigma$  contains all open sets and
is closed under countable unions.  We want to show that $\Sigma$
is also closed under complements.

For a set $G \in \cal G$ let 
$$G'=\{\<x,y\> : y \text{ is an interior point of } 2^\omega \setminus
(G)_x\}.$$

Note that $(2^\omega \times 2^\omega) \setminus (G \cup G')$ is a set
whose vertical sections are  closed and nowhere dense.
It follows that in order to show  that $\Sigma $ is closed under
complements it is enough to check that $G'$ is a Borel set.

Let $\{U_n :n \in \omega\}$ be a recursive enumeration of a countable
base for the family 
of open subsets of 
$2^\omega$. 

Note that the following are equivalent:
\begin{enumerate}
\item $\<x,y\> \in G'$,
\item $ \exists n \ \lft2(y \in U_n \ \& \ \forall z \
\lft1(z \not\in U_n \ \vee\  \<x,z\> \not\in G\rgt1)\rgt2)$
(${\boldsymbol \Pi}^1_1$),
\item $\exists n \ \lft2(y \in U_n \ \& \ \forall m \ \lft1( U_n \cap U_m =
\emptyset \ \vee\  \exists z \ (z \in U_m \ \<x,z\> \not\in
G\rgt1)\rgt1)\rgt2)$ (${\boldsymbol \Sigma}^1_1$).
\end{enumerate}

That shows that $G'$ has a $\boldsymbol \Delta^1_1$ definition which
means that it is a Borel set.

\vspace{0.1in}

(2) Let $\Sigma$ be the collection of all sets $H$ satisfying the
 conclusion of the theorem.
We will show that $\Sigma$ is a $\sigma$-algebra containing all open
sets. In particular, all  Borel sets are in $\Sigma$.
It is enough to show that 
\begin{enumerate}
  \item[(i)] unions of finitely many rectangles are in $\Sigma$,
  \item[(ii)] if $A_0 \subseteq A_1 \subseteq A_2 \subseteq \cdots$ are in
    $\Sigma$ then $\bigcup_{n \in \omega} A_n \in \Sigma$,
  \item[(iii)] if $A_0 \supseteq A_1 \supseteq A_2 \supseteq \cdots$ are in
    $\Sigma$ then $\bigcap_{n \in \omega} A_n \in \Sigma$.
\end{enumerate}

Condition (i) is clear. To show (ii) fix $\varepsilon>0$ and let $B_n$
witness that $A_n \in \Sigma$ with $\varepsilon_n=\varepsilon
2^{-n}$. Let $B=\bigcup_{n \in \omega} B_n$.

\vspace{0.1in}

(iii) Fix $\varepsilon>0$ and let $B_n$
witness that $A_n \in \Sigma$ with $\varepsilon/2$. 
For $n \in \omega$ let $$Z_n=\left\{x: \mu\lft1((A_n \setminus A)_x\rgt1) <
\varepsilon/2\right\}.$$ 
Put $B=\bigcup_{n \in \omega} \lft1((Z_n \times 2^\omega) \cap
B_n\rgt1)$.~$\QED$

\section{$SR^{\M}$ sets}
In this section we will characterize 
$SR^{\M}$ sets.

We will need the following characterization of meager sets.

\begin{theorem}[\cite{BJS} Prop. 9]\label{representmeager}
  For every meager set $F \subseteq 2^\omega$ there exists $x_F \in
  2^\omega$ and a strictly increasing function $f_F \in \omega^\omega$
  such that 
$$F \subseteq  B(f_F,x_F)=\{x \in 2^\omega : \forall^\infty n \ \exists j \in
\lft1[f_F(n), f_F(n+1)\rgt1) \ x(j) \neq x_F(j)\}.$$

Moreover,
$B(f,x) \subseteq B(g,y)$ if and only if
\begin{multline*}
\forall^\infty n \ \exists k \ \lft2(g(n)\leq f(k)<f(k+1)\leq g(n+1) \
\& \\
x\rest\lft1[f(k),f(k+1)\rgt1) = y\rest\lft1[f(k),f(k+1)\rgt1)\rgt2).~\QED
\end{multline*}
\end{theorem}

\begin{theorem}\label{charmeaadd}
\begin{enumerate}
\item $\add(\M)$  is the least size of any family $F \subseteq \omega^\omega$
such that  
there are no $r,h \in \omega^\omega$ such that 
$$\forall f \in F \ \forall^\infty n \
\exists k \in  \lft1[r(n), r(n+1)\rgt1) \ f(k) = h(k),$$
\item $X$ is $\M$-additive iff 
for every increasing function $f \in \omega^\omega$ there exists $g
\in \omega^\omega$ and $y \in 2^\omega$ such that 
\begin{multline*}
\forall x \in X \ \forall^\infty n \ \exists k \ \lft2(g(n)\leq
f(k)<f(k+1)\leq g(n+1) \ 
\& \\
x\rest\lft1[f(k),f(k+1)\rgt1) = y\rest\lft1[f(k),f(k+1)\rgt1)\rgt2),
\end{multline*}
 \item $\unif(\M^\star)$  is the least size of a bounded  family $F \subseteq
 \omega^\omega$ 
such that  
there are no $r,h \in \omega^\omega$ such that 
$$\forall f \in F \ \forall^\infty n \
\exists k \in  \lft1[r(n), r(n+1)\rgt1) \ f(k) = h(k).$$
\item $\add(\M) = \min\{\unif(\M^\star), {\frak b}\}.$
\end{enumerate}
\end{theorem}
\Proof 
(1) Follows readily from the Miller-Truss result that $\add(\M)=\min\{\cov(\M), \frak
b\}$ and the fact that $\cov(\M)$ is the least cardinal of any set $F
\subseteq \omega^\omega$ such that there is no $h \in \omega^\omega$
such that 
$$\forall f \in F \ \exists^\infty n \ h(n)=f(n).$$
 See \cite{Mil81Som} and  \cite{Bar87Com}.

(2) Follows from \ref{representmeager} and the fact that $B(f,x)+z=B(f, x+z)$.

(3) Similar to (2).

(4) It was proved by Pawlikowski in \cite{Paw85Pow}. It follows from
(1) and (3) and the fact that $\add(\M) \leq \frak b$.~$\QED$

We will now characterize $SR^{\M}$ sets.
\begin{theorem}\label{first}
  A set $X \subseteq 2^\omega$ is an $SR^{\M}$ set iff every
   Borel image of $X$ is $\M$-additive and every Borel image of $X$
   into $\omega^\omega$ is bounded.
\end{theorem}
\Proof By \ref{1.2},  $SR^{\M}$ sets are
$\M$-additive. Thus by \ref{clear}, all Borel images of
$SR^{\M}$ sets  are $\M$-additive. 

To show the second part let $F$ be a homeomorphism (Borel isomorphism
is enough)
between $2^\omega
\setminus \rationals$ and the set of increasing functions in $\omega^\omega$.
Let $H = \{\<x,y\> : y \in B(F(x),x)\}$. Suppose that $f$ is a Borel
mapping of $X$ into $\omega^\omega$ such that the family $f(X)$ is
unbounded. Without loss of generality we can assume that $f(X)$
consists of increasing functions. Let $Y=F^{-1}(f(X))$. By
\ref{clear}, $Y$ is an $SR^{\M}$ set and by
\ref{representmeager}, 
$$\bigcup_{x \in Y} (H)_x \not \in \M . $$
This contradiction finishes the proof of one implication.

\vspace{0.1in}

($\leftarrow$)
Suppose that $H \subseteq 2^\omega \times 2^\omega$ is a Borel $\M$-set.

Thus, by \ref{fremlin}(1) we can assume that $H=\bigcup_{n \in
  \omega} G_n$, where each set $G_n$ has closed nowhere dense
sections.
Moreover, by combining \ref{representmeager} and \ref{fremlin}(1), we can find
a Borel mapping $F$ such that for every $z \in 2^\omega$,
$F(z)=\<f_z, x_z\> \in \omega^\omega \times 2^\omega$ is such that 
$(H)_z \subseteq B(f_z,x_z)$. 
By the assumption, the family $\{f_z : z \in X\}$ is bounded.

Without loss of generality we can assume that $f_z = f$ for $z \in X$.
By \ref{charmeaadd}(2), there exists a function $g \in \omega^\omega$
and real $y \in 2^\omega$
such that 
\begin{multline*}
\forall z \in X \ \forall^\infty n \ \exists k \ \lft2(g(n)\leq f(k)<f(k+1)\leq g(n+1) \
\& \\
x_z\rest\lft1[f(k),f(k+1)\rgt1) = y\rest\lft1[f(k),f(k+1)\rgt1)\rgt2).
\end{multline*}

It is clear that 
$$\bigcup_{z \in X} (H)_z \subseteq B(g,y) \in \M, $$
which finishes the proof.~$\QED$

\section{$SR^{\N}$ sets}
In this section we characterize $SR^{\N}$ sets. 

We start with the following well-known fact:
\begin{theorem}[\cite{Bar84Add}, \cite{Paw85Pow}]
 $\add(\N)$ is the least size of any  family $F
  \subseteq \omega^\omega$ such that there is no function $S : \omega
\longrightarrow [\omega]^{<\omega}$ with $|S(n)| \leq n$ for all $n$,
such that 
$$\forall f \in F \ \forall^\infty n \ f(n) \in S(n).$$
  $\unif(\N^\star)$ is the least size of any bounded family $F
  \subseteq \omega^\omega$ such that there is no function $S : \omega
\longrightarrow [\omega]^{<\omega}$ with $|S(n)| \leq n$ for all $n$,
such that 
$$\forall f \in F \ \forall^\infty n \ f(n) \in S(n).$$

In particular, $\add(\N) = \min\{\unif(\N^\star), {\frak b}\}$.~$\QED$
\end{theorem}
Note that the characterization of $\unif(\N^\star)$ above is an easy
corollary of the equivalence between (1) and (3) in the theorem below.

We get a  similar characterisation of $\N$-additive sets but the
proof is much harder.
We will present it here for completeness.

\begin{theorem}[Shelah \cite{Sh445}]\label{charnuladd}
Let $X \subseteq 2^\omega$. The following conditions are equivalent:
\begin{enumerate}
\item $X$ is $\N$-additive,
\item for every increasing function $f \in \omega^\omega$ there exists
  a tree $T \subseteq 2^{<\omega}$ such that for all $n$, $|T \rest n|\leq
  f(n)$ and for every $x \in X$ there 
  exists $n$ such that $x \in [T^{\<n\>}]$,
\item for every increasing function $f \in \omega^\omega$ there exists
  a sequence $\{I_n : n \in \omega\}$ such that 
  \begin{enumerate}
  \item for all $n$, $I_n \subseteq 2^{[f(n), f(n+1))}$,
  \item for all $n$, $|I_n| \leq n$,
  \item $\forall x \in X \ \forall^\infty n \ x \rest \lft1[f(n),
    f(n+1)\rgt1) \in I_n$.
  \end{enumerate}
\item there exists a function $g \in \omega^\omega$ such that for
  every increasing function $f \in \omega^\omega$ there exists 
  a sequence $\{I_n : n \in \omega\}$ such that 
  \begin{enumerate}
  \item for all $n$, $I_n \subseteq 2^{[f(n), f(n+1))}$,
  \item for all $n$, $|I_n| \leq g(n)$,
  \item $\forall x \in X \ \forall^\infty n \ x \rest \lft1[f(n),
    f(n+1)\rgt1) \in I_n$.
  \end{enumerate}
\end{enumerate}
\end{theorem}
\Proof $(1) \rightarrow (4)$  The function $g$ that we are looking for
will be $g(n)=n 2^n$.

Let $f \in
\omega^\omega$ be an increasing function such that for all 
$n$,
$$f(n+1) \geq 2^{f(0)+ \cdots + f(n) + n}.$$
We will start with a construction of a measure zero set having strong
combinatorial properties. \label{shelahset}
\begin{lemma}\label{trivlem}
  For every $n$ there exists a family $\{A_s : s \in 2^{f(n)}\}$ such
  that 
  \begin{enumerate}
  \item $A_s \subseteq 2^{[f(n), f(n+1))}$ for $s \in 2^{f(n)}$,
  \item $\mu\lft1([A_s]\rgt1) =  1-2^{-n}$ for $s \in 2^{f(n)}$,
  \item sets $\{x_s+[A_s] : s \in 2^{f(n)}\}$ are
probabilistically independent for every family $\{x_s : s \in
2^{f(n)}\} \subseteq 2^\omega$.
  \end{enumerate}
\end{lemma}
\Proof Such a family may be constructed in many different ways.  Below
is one such construction.

Fix a family of sets $\{I_s : s \in 2^{f(n)}\}$ such that 
\begin{enumerate}
\item $I_s \subseteq \lft1[f(n), f(n+1)\rgt1)$ for all $s$,
\item $I_s \cap I_t = \emptyset$ for $s \neq t$,
\item $|I_s|=n$ for all $s$.
\end{enumerate}
For $s \in 2^{f(n)}$ let 
$$A_s =\left\{t \in 2^{[f(n),f(n+1))} : \exists k \in I_s \
t(k)=1\right\}.~\QED$$ 

Define a sequence of trees $\{T_m: m \in \omega  \}$  using the
following condition: 

$$T_m \rest f(m) = 2^{f(m)} \text{ and for } n \geq m
\text{ if } s \in T_m \rest {f(n)} \text{ then }
\suc_{T_m,f(n+1)}(s)=A_s .$$ 

Note that
$$\mu\lft1([T_m]\rgt1) = \prod_{n=m}^\infty \left(1-2^{-n}\right)
\stackrel{m \rightarrow \infty}{\longrightarrow} 1 .$$
In particular, the set $H = 2^\omega \setminus \bigcup_{m \in \omega}
[T_m]$ has measure zero.
Since $X$ is a null additive the set $X+H$ has measure zero. Passing to
complements we conclude that the set 
$$\bigcap_{x \in X} \left(x + \bigcup_{m \in \omega}
[T_m]\right)$$ has positive measure. 

Let $T^\star \subseteq
2^{<\omega}$ be a tree such that $\mu\lft1([T^\star_s]\rgt1)>0$ for
all $s \in T^\star$ and 
$$[T^\star] \subseteq \bigcap_{x \in X} \left(x
+ \bigcup_{m \in \omega} [T_m]\right).$$

Note that if $[T^\star] \subseteq 
\bigcup_{m \in \omega} \left[T_m+x\right]$ then there
is $s \in T^\star$ and $m \in \omega$  such that $T^\star_s \subseteq
T_m+x$.  

For $s \in T^\star$ and $m \in \omega$ define 
$$X_{s,m} = \{x \in X : T^\star_s \subseteq T_m+x \}.$$
By the above remark $X = \bigcup_{s \in T^\star, m \in \omega}
X_{s,m}$.

We will first show that each set $X_{s,m}$ satisfies the conditions of
(4).

Fix $\bar{s} \in T^\star, \bar{m} \in \omega$ and let $n  \geq \bar{m},
|\bar{s}|$. By extending, if necessary, we can assume that
$\mu([T^\star_{\bar{s}}])> \frac{1}{2} \mu([\bar{s}])$.

We will estimate the size of the set $\{x \rest
f(n) : x \in X_{\bar{s},\bar{m}}\}$.
Consider the finite  tree 
$$\widetilde{T}=\bigcap_{x \in X_{\bar{s},\bar{m}}}
\lft2(\lft1(T_{\bar{m}} \rest 
f(n+1)\rgt1) + x \rest f(n+1)\rgt2).$$
Observe that $\widetilde{T}$ has height $f(n+1)$ and contains
$T^\star_{\bar{s}} \rest 
f(n+1)$. 
Let $\bar{t}$ be any element of $\widetilde{T} \rest f(n)$. 
For every element $x  \in X_{\bar{s},\bar{m}}$ there is a
$t_x \in 
T_{\bar{m}} \rest f(n)$ such that $x \rest f(n) +t_x =\bar{t}$.
Since $t_x$ depends only on $x \rest f(n)$ the sets $\{x \rest f(n) :
x \in X_{\bar{s},\bar{m}} \}$ and $\{t_x : x \in
X_{\bar{s},\bar{m}}\}$ have the same size. 
It follows that 
$$\suc_{\tilde{T}, f(n+1)}(\bar{t}) = \bigcap \left\{
\suc_{T_{\bar{m}}, f(n+1)}(t_x) + x \rest f(n+1) : x \in
X_{\bar{s},\bar{m}}\right\}.$$ 
The sequences $t_x$ are all different and therefore sets
$\suc_{T_{\bar{m}},f(n+1)}(t_x)$ represent probabilistically
independent sets.  Thus their translations are also independent thus
computing measures (relativised to $[\bar{t}]$) and using the fact
that the sets on the right hand 
side are independent we get
$$\frac12  <
\mu\lft1([T^\star_{\bar{t}}]\rgt1) \leq 
\mu\lft2(\lft1[\suc_{\tilde{T}, f(n+1)}(\bar{t})\rgt1]\rgt2) =
\left(1-\frac{1}{2^n}\right)^{|\{x \rest f(n) : x \in
  X_{\bar{s},\bar{m}}\}|} .$$ 
In particular we get that for all $n \geq |\bar{s}|,\bar{m}$,
$$|\{x \rest f(n) : x \in
X_{\bar{s},\bar{m}}\}| \leq 2^n .$$ 
Let $\{X^k : k \in \omega\}$ be enumeration of the family $\{X_{s,m} :
s \in T^\star, m \in \omega\}$. For each $k,m \in \omega$ let 
$I^k_m = \{x \rest \lft1[f(m), f(m+1)\rgt1) : x \in X^k\}$.
By the above argument
$$\forall k \ \forall^\infty m \ |I^k_m| \leq 2^{m+1} .$$
Let $g(n)=n\cdot 2^{n+1}$ for all $n$. Define for $n \in \omega$,
$$I_n = \bigcup \left\{I^k_n : k \leq n, |I_n^k| \leq
2^{n+1}\right\}.$$
It is clear that the sequence $\{I_n:n \in \omega\}$ has the 
required properties.

Note that we have proved the following:
\begin{lemma}
Suppose that $f \in \omega^\omega$ is strictly increasing and that for
each $n$ we have a disjoint family $\{I_s: s \in 2^{f(n)}\}$ of subsets
of $[f(n), f(n+1))$ of size $n$. 
Let 
$$H=\{x \in 2^\omega: \exists^\infty n \ \forall j \in I_{x \rest
  f(n)} \ x(j)=0\}.$$
Then $H$ has measure zero and if $X$ is any set such that $X+H$ has
measure zero then there are sets $J_n \subseteq 2^{f(n)}$ such that
$|J_n| \leq n2^n$ and
$$\forall x \in X \ \forall^\infty n \ x \rest f(n) \in J_n.~\QED$$
\end{lemma}

\vspace{0.1in}

$(4) \rightarrow (3)$ Suppose that the
function $f \in \omega^\omega$ is given.
Apply  (4) to the function $f'(n)=f\lft1(g(n)\rgt1)$ for $n \in
\omega$. 

\vspace{0.1in}

Implication $(3) \rightarrow (2)$ is very easy.

\vspace{0.1in}

$(2) \rightarrow (1)$ Let $H \subseteq 2^\omega $ be a measure zero
set. 
The following lemma is well-known:
\begin{lemma}\label{oxtoby} 
There exists a sequence $\< F_{n} : n  \in
\omega  \>$  such
that  $F_{n}  \subseteq  2^{n}$ for $n  \in
\omega$,  $\sum_{n=1}^{\infty}   |F_{n}|\cdot2^{-n} <
\infty$   and  $H  \subseteq  \left\{x  \in     2^{\omega}  :
\exists^{\infty} n \ x
\rest n  \in  F_{n} \right\}$.~$\QED$
\end{lemma}
By \ref{oxtoby}, there exists a sequence $\< F_{n} : n  \in
\omega  \>$  such
that  $F_{n}  \subseteq  2^{n}$ for $n  \in
\omega$,  $\sum_{n=1}^{\infty}   |F_{n}|\cdot2^{-n} <
\infty$   and  $H  \subseteq  \bigcap_{m \in \omega} \bigcup_{k>m} [F_k]$.

Let $f \in \omega^\omega$ be an nondecreasing function with $\lim_{n
\rightarrow \infty } f(n)=\infty$ such that 
$$\sum_{n=1}^{\infty}   \frac{f(n)\cdot |F_{n}|}{2^{n}} < \infty.$$

Let $T$ be a tree from \ref{charnuladd}(2) for this function. Note
that $|T^{\<n\>} 
\rest k| \leq 2^n f(k)$ for all $n, k \in \omega$. 

We have
$X+H \subseteq \bigcup_{n \in \omega} [T^{\<n\>}] + H$.
For every $n$,
$$[T^{\<n\>}] + H \subseteq \bigcap_{m \in \omega} \bigcup_{k > m}
[T^{\<n\>}] + [F_k].$$
By the choice of $f$ the measure of  $[T^{\<n\>}] + [F_k]$ is bounded by
$2^n \cdot f(k)\cdot  |F_k|\cdot 2^{-k}$. 
Thus
$$\mu\left([T^{\<n\>}] + H\right) \leq 2^n\cdot \sum_{k=m}^\infty
\frac{f(k)\cdot 
  |F_k|}{2^{k}} \stackrel{m \rightarrow \infty}{\longrightarrow } 0
.$$
Since $n$ is arbitrary we conclude that $X+H$ has measure zero which
finishes the proof.~$\QED$

As a consequence we have:
\begin{theorem}[Shelah \cite{Sh445}]
  $\N^\star \subseteq \M^\star$. 
\end{theorem}
\Proof Follows immediately from \ref{charnuladd} and
\ref{charmeaadd}.~$\QED$

\begin{theorem}\label{charstrmead}
  A set $X \subseteq 2^\omega$ is a $SR^{\N}$ set iff every
   Borel image of $X$ is $\N$-additive and every Borel image of $X$
   into $\omega^\omega$ is bounded.
\end{theorem}
\Proof As in \ref{first}, we show that the Borel images of
$SR^{\N}$ sets are $\N$-additive.
To show the second part let $F$ be a homeomorphism between $2^\omega
\setminus \rationals$ and the set of increasing functions in
$\omega^\omega$.
For an increasing function $f \in \omega^\omega$ let 
$$B(f)=\{x \in 2^\omega : \exists^\infty n \ \forall i \leq n \
x(f(n)+i)=0\}.$$ 
Clearly, $B(f)$ is a measure zero set.
Let $H = \{\<x,y\> : y \in B(F(x))\}$. Suppose that $f$ is a Borel
mapping of $X$ into $\omega^\omega$ such that the family $f(X)$ is
unbounded. Without loss of generality we can assume that $f(X)$
consists of increasing functions. Let $Y=F^{-1}(f(X))$. By
\ref{clear}, $Y$ is an $SR^{\N}$ set. By a result of Miller
(\cite{Mil84Add}, lemma 5),
$$\bigcup_{x \in Y} (H)_x \not \in \N.$$

\vspace{0.1in}

To show the other implication  we need the following fact.
\begin{theorem}[Rec{\l}aw \cite{RecOpen}]\label{reclaw}
Let $X \subseteq 2^\omega$ be a set such that for every Borel function
$x \leadsto f^x \in \omega^\omega $ there exists a function $S :
\omega \longrightarrow [\omega]^{<\omega}$ such that $|S(n)|\leq n$
for all $n$ 
and
$$\forall x \in X \ \forall^\infty n \ f^x(n) \in S(n) .$$
Then $X$ is a $SR^{\N}$ set.~$\QED$
\end{theorem}

Now we are ready to finish the proof of the theorem \ref{charstrmead}.
Suppose that $F$ is a Borel  mapping of $X$ into $\omega^\omega$. By the
assumption the set $F(X)$ is bounded.
Let $f$ be a function such that $F(X)$ is bounded by the function
$2^{f(n+1)-f(n)}$. Identify $2^{f(n+1)-f(n)}$ with $2^{[f(n),f(n+1))}$
for all $n$. In this way we can identify $F(X)$ with a subset of
$2^\omega$. Part (3) of \ref{charnuladd} and \ref{reclaw} conclude the
proof.~$\QED$ 

Note that the assumptions that Borel images of $X$ into $\omega^\omega$
are bounded in \ref{charmeaadd} and \ref{charstrmead} are
necessary. It follows immediately from the following theorem:
\begin{theorem}[Rec{\l}aw \cite{Rec89Sma}]\label{reclaw1}
  Assume Martin's Axiom. Then the real line is a Borel image of some
  $\N$-additive set $X$.~$\QED$
\end{theorem}

In particular, the set $X$ from \ref{reclaw1}, is $\N$- and
$\M$-additive but is neither a $SR^{\N}$- nor $SR^{\M}$ set.

It is also consistent (see \cite{Paw85Pow}) that there exists a set $X
\not \in SR^{\N}$ such that all Borel images of $X$ are null-additive.

\section{$R^{\M}$ sets}

In this section we will study $R^{\M}$ sets. All the results of this
section can be found in \cite{BarJud93Cov}. They were also proved by
Pawlikowski and Rec{\l}aw in \cite{PawRecPar} and \cite{PawPropC} and
Rec{\l}aw in \cite{RecOpen}.
\begin{definition}
  A set $X \subseteq 2^\omega$ has Rothberger's  property (is a $C''$
  set) if for every sequence of open covers of $X$, $\{{\cal G}_n :n
  \in \omega\}$ there exists a sequence $\{U_n : n \in \omega\}$ with
  $U_n \in {\cal G}_n$ such that $X \subseteq \bigcup_{n \in \omega} U_n$.
\end{definition}

Note that Rothberger's property is the topological version of strong
measure zero. We have the following:
\begin{theorem}[Fremlin, Miller \cite{FreMil88Som}]
  The following are equivalent:
  \begin{enumerate}
  \item $X \subseteq 2^\omega $ has Rothberger's property
\item $X$ has strong measure zero with respect to every metric which
  gives $X$ the same topology.~$\QED$
  \end{enumerate}
\end{theorem}

Let $C''$ be the collection of subsets of $2^\omega$ which have
Rothberger's property. It is easy to see that $C''$ is a $\sigma$-ideal.

\begin{theorem}[Bartoszy{\'{n}}ski, Judah, Pawlikowski, Rec{\l}aw
  \cite{BarJud93Cov}, \cite{RecOpen}, \cite{PawRecPar}]\label{main}
  The following conditions are equivalent:
  \begin{enumerate}
  \item $X$ is an $R^{\M}$ set,
  \item for every Borel function $x \leadsto f^x \in \omega^\omega$
    there exists a function $g \in 
  \omega^\omega$ such that 
$$\forall  x \in X \ \exists^\infty n \
f^x(n)=g(n).$$
\item for every Borel function 
$x \leadsto \<Y^x,f^x\> \in [\omega]^\omega \times \omega^\omega$
there exists $g \in  
\omega^\omega$ such that  
$$\forall x \in X \ \exists^\infty n \in Y^x \ f^x(n)=g(n).$$
\item Every Borel image of $X$ has Rothberger's property.~$\QED$
  \end{enumerate}
\end{theorem}

As a consequence we get:
\begin{theorem}[Bartoszy{\'{n}}ski, Judah \cite{BarJud93Cov}]
  $\cf\lft1(\cov(\M)\rgt1) \geq \add(R^{\M}) \geq \add(C'') 
  \geq \add(\N)$.~$\QED$ 
\end{theorem}

\section{$R^{\N}$ sets }
In this section we will study $R^{\N}$ sets. 
Most of the results of this section were obtained independently (and
earlier) by Pawlikowski.

Let us start with the
following definition.
\begin{definition}
  A set $G \subseteq 2^\omega$ is called small if
there exists a sequence of disjoint intervals $\{I_n: n \in \omega\}$
and a sequence $\{J_n : n \in \omega\}$ such that for all $n \in \omega$,
\begin{enumerate}
\item $J_n \subseteq 2^{I_n}$,
\item $|J_n|\cdot 2^{-|I_n|} \leq 2^{-2n}$,
\item $G \subseteq \{x \in 2^\omega : \exists^\infty n \ x \rest I_n
  \in J_n\}$.
\end{enumerate}
We denote the set $\{x \in 2^\omega : \exists^\infty n \ x \rest I_n
  \in J_n\}$ by $(I_n,J_n)_{n=0}^\infty$.
 
For a function $ f \in \omega^\omega$ such that $f(n)>0$ for all $n$
and a series $\sum_{n=0}^\infty \varepsilon_n<\infty$ let 
$${\cal X}_f = \prod_{n \in \omega} f(n)$$
 and 
$$\Sigma_f =\left\{S \in \lft1([\omega]^{<\omega}\rgt1)^\omega :
\frac{|S(n)|}{f(n)} < \varepsilon_n\right\}.$$
\end{definition}
For the rest of this section we will assume that
$\varepsilon_n=2^{-n}$ for all $n$. 

Note that by setting $f(n) = 2^{|I_n|}$ and $S(n)=J_n$ we can identify
the set 
$$\{x \in 2^\omega : \exists^\infty n \ x \rest I_n
  \in J_n\}$$
 with the set $\{z \in {\cal X}_f : \exists^\infty n \
  z(n) \in S(n)\}$.

\begin{definition}
Let ${\cal H}$ be the ideal of all sets $X \subseteq 2^\omega$ such
that every Borel image of $X$ into $\omega^\omega$ is bounded. 
\end{definition}

By replacing ``Borel'' by ``continuous'' we get a weaker property
introduced by Hurewicz as $E^{\star\star}$ (\cite{HuSeq27})

We will start with the following forcing characterization of $R^{\N}$
sets.

Recall that for a model $M$, $\R(M)$ denotes the set of random reals
over $M$.

\begin{theorem}
  $X$ is an $R^{\N}$-set iff for
  every countable elementary submodel model $M \prec H(\lambda)$ 
there exists a real $z$ such
  that $z$ is random over $M[x]$ for $x \in X$.
\end{theorem}
\Proof  Recall that $H(\lambda)$ is the collection of sets
hereditarily of size $<\lambda$. If $\lambda$ is a regular,
uncountable cardinal then $H(\lambda)$ is a model for a large fragment
of $\ZFCa$. Moreover, in this context, $M[x]$ is a ``closure'' of $M
\cup \{x\}$.

$(\leftarrow)$ Let $H$ be an $\N$-set. Choose $M$
containing the code of $H$ and let $z$ be as above. Then $z \not \in
\bigcup_{x \in X} (H)_x$.

\vspace{0.1in}

$(\rightarrow)$ Let $M \thinks \ZFCs$ be a countable model. 
Consider the set
$$\bar{H} = \left\{\<x,y\> : y \not \in \R\lft1(M[x]\rgt1)\right\}.$$
$M$ is countable and can be coded as a real. Note that $y \not \in
\R\lft1(M[x]\rgt1)$  is equivalent to the statement ``there exists a
code for a null set $F \in M[x]$ such that $y \in F$''. Since $M$ is
countable this statement is arithmetical (in $M$). It follows that
$\bar{H}$  is a Borel set.  Let $z$ be such that $z \not \in
(\bar{H})_x$ for $x \in X$.  
Clearly $z$ is the real we are looking for.~$\QED$

\begin{theorem}\label{oldstuff}
  Suppose that $X \in {\cal H}$. The following are equivalent:
  \begin{enumerate}
  \item X is an $R^{\N}$ set,
\item for every function $f \in \omega^\omega$ and every Borel mapping
  $x \leadsto S^x \in \Sigma_f$ there exists a function $g \in
  {\cal X}_f$ 
  such that 
$$\forall x \in X \ \forall^\infty n \ g(n) \not \in S^x(n).$$
  \end{enumerate}
\end{theorem}
\Proof Implication $(1) \rightarrow (2) $ is very easy because
$H=\{(x,g): \exists^\infty n \ g(n) \in S^x(n)\}$ is a Borel $\N$-set
in $X \times {\cal X}_f$.

\vspace{0.1in}

$(2) \rightarrow (1)$
The proof is based on techniques from \cite{Bar88Cov}. 

Suppose that $H
\subseteq 2^\omega\times 2^\omega$ is a Borel $\N$-set.

\begin{lemma}
  There exist increasing interleaved sequences $\{n_k,m_k : k \in \omega\}$ and
  families $\{J^x_k, \widetilde{J}^x_k : k \in \omega\}$ such that for
  all $k \in \omega$,
  \begin{enumerate}
  \item $J^x_k \subseteq 2^{[n_k, n_{k+1})}$, $\widetilde{J}^x_k
    \subseteq 2^{[m_k, m_{k+1})}$,
\item mappings $x \leadsto \{J^x_k : k \in \omega\}$ and $x \leadsto
  \{\widetilde{J}^x_k : k \in \omega\}$ are Borel,
\item $|J^x_k|\cdot 2^{n_k - n_{k+1}} \leq 2^{-2k}$,
    $|\widetilde{J}^x_k|\cdot 2^{m_k - m_{k+1}} \leq 2^{-2k}$ for all
      $x \in X$,
    \item $(H)_x \subseteq \lft1([n_k,n_{k+1}),J^x_k\rgt1)_{k=0}^\infty
      \cup \lft1([m_k,m_{k+1}),\widetilde{J}^x_k\rgt1)_{k=0}^\infty$.
  \end{enumerate}
\end{lemma}
\Proof
By \ref{fremlin}(2), we can find a sequence of Borel  sets $\{B_n : n \in
\omega\}$ such that $(H)_x \subseteq (\bigcap_{n \in \omega} B_n)_x$
for $x \in X$ and 
$$\forall x \in X \ \forall^\infty n \ \mu\lft1((B_n)_x\rgt1) <
\frac{1}{2^n}.$$ 
Thus we can work with $B_n$'s rather than with $H$.

The rest of the proof is the repetition of the proof of
 theorem 2.2 in \cite{Bar88Cov} (using the fact that
$X \in {\cal H}$).~$\QED$

As in \cite{Bar88Cov}, define for $k \in \omega$ and $x \in X$ sets
$$S^x_k = \left\{s \in 2^{[n_k,m_k)} : s \text{ has at least
  $2^{n_{k+1}-m_k-k}$ extensions in } J^x_k\right\}$$ and for $k>0$
$$\widetilde{S}^x_k = \left\{s \in 2^{[n_k,m_k)} : s \text{ has at least
  $2^{n_{k}-m_{k-1}-k}$ extensions in } \widetilde{J}^x_{k-1}\right\}.$$ 

The mapping $x \leadsto \{S^x_k \cup \widetilde{S}^x_k : k \in \omega\}$
is Borel. Moreover,  $\{S^x_k \cup \widetilde{S}^x_k : k \in \omega\}
\in \Sigma_f$, where $f(k)=2^{m_k-n_k}$. Therefore, by the assumption
about $X$ there exists a real $z$ such that 
$$\forall x \in X \ \forall^\infty k \ z \rest [n_k,m_k) \not \in
S^x_k \cup \widetilde{S}^x_k.$$
Define for $x \in X$ and $k \in \omega$,
$$T^x_k = \left\{s \in 2^{[m_k,n_{k+1})} : z \rest [n_k,m_k)^\frown s \in
  J^x_k \text{ or } s^\frown z \rest [n_{k+1},m_{k+1}) \in
  \widetilde{J}^x_{k+1}\right\}.$$ 
As before we easily check that the assumption about $X$ yields that
there exists a real $y$ such that 
$$\forall x \in X \ \forall^\infty k \ y \rest [m_k,n_{k+1}) \not \in
T^x_k.$$

Now define $\widetilde{z} \in 2^\omega$ as
$$\widetilde{z}(n) = \left\{\begin{array}{ll}
z(n) & \text{if $n \in [n_k,m_k)$ for some $k$}\\
y(n) & \text{if $n \in [m_k,n_{k+1})$ for some $k$}
\end{array}\right. .$$
As in the proof of 2.2 in \cite{Bar88Cov} we check that 
$$\widetilde{z} \not \in \bigcup_{x \in X} (H)_x .~\QED$$

\begin{theorem}[Pawlikowski \cite{Paw92Sie}]
Suppose that $X \in {\cal H}$. Then $X \in R^{\N}$ iff for every
measure zero set $E \in \N$, $E \cap X \in R^{\N}$.
\end{theorem}
\Proof
The implication $(\rightarrow) $ is obvious.

\vspace{0.1in}

To show the other implication consider the space ${\cal X}_f=\prod_{n
  \in \omega} f(n)$ equipped with
its standard product measure. 

We will work in the space $2^\omega \times {\cal X}_f$.

By \ref{oldstuff} it is enough to show that 
for every function $f \in \omega^\omega$ and every Borel mapping
  $x \leadsto S^x \in \Sigma_f$ there exists a function $g \in {\cal X}_f$
  such that 
$$\forall x \in X \ \forall^\infty n \ g(n) \not \in S^x(n).$$

Suppose that a function $f$ and a Borel mapping as above are given.
Let $H \subseteq 2^\omega \times {\cal X}_f$ be the Borel set such that 
$(H)_x = \{h \in {\cal X}_f : \exists^\infty n \ h(n) \in S^x(n)\}$.

Construct a sequence of elements of ${\cal X}_f$, $\{g_n : n \in
\omega\}$  and a sequence of measure zero sets $\{E_n : n \in
\omega\}$ such that 
\begin{enumerate}
\item $E_n = \{x\in X : g_n \in (H)_x\}$ for all $n$,
\item $g_{n+1} \in {\cal X}_f \setminus \bigcup \{(H)_x : x \in E_0 \cup
  \cdots \cup  E_n\}$.
\end{enumerate}
The existence of these sequences follows from Fubini's theorem and the
assumption about $X$. Note that the set $\bigcup \{(H)_x : x \in E_0 \cup
  \cdots \cup  E_n\}$ does not have full measure.

For $x \in X$ define $h^x \in \omega^\omega$ as follows:
$$h^x(n)=\min\{ k : \forall j > k \ g_n(j) \not \in S^x(j)\} \text{ for
  } n \in \omega.$$
Note that $h^x(n)$ is defined for all except finitely many values of
$n$.

Since the mapping $x \leadsto h^x$ is Borel  we can find an increasing
function $h \in \omega^\omega$ such that $h^x \leq^\star  h$ for all $x \in
X$.
Let 
$$g(k) = g_n(k) \text{ for } h(n)<k \leq h(n+1).$$
It is easy to see that 
$$\forall x \in X \ \forall^\infty n \ g(n) \not \in S^x(n)$$
which finishes the proof.~$\QED$

As a corollary we get:
\begin{theorem}
  Suppose that $X \cap E \in SR^{\N}$ for every $E \in \N$. Then $X
  \in R^{\N}$.
\end{theorem}
\Proof
It is enough to show that $X \in {\cal H}$. Suppose that $F: 2^\omega
\longrightarrow \omega^\omega$ is a Borel mapping. We can find a
sequence of compact sets $\{A_n : n \in \omega\}$ such that 
\begin{enumerate}
\item $F \rest A_n $ is continuous,
\item $\bigcup_{n \in \omega} A_n$ has full measure. 
\end{enumerate}
Let $E = 2^\omega \setminus \bigcup_{n \in \omega} A_n$. Since $E \cap X \in
SR^{\N}$, it follows from \ref{charstrmead} that $F(E\cap X)$ is
bounded. On the other hand, $F(\bigcup_{n \in \omega}A_n)$ is also
bounded, which finishes the proof.~$\QED$

A set $X \subseteq 2^\omega$ is a Sierpi{\'{n}}ski set if $X \cap H$ is
countable for every measure zero set $H \subseteq 2^\omega$.

\begin{theorem}[Pawlikowski \cite{Paw92Sie}]
Every Sierpi{\'{n}}ski set is strongly meager.
\end{theorem}
\Proof All countable sets are in $SR^{\N}$ and all $R^{\N}$ sets are
strongly meager.~$\QED$

We do not know if the  $R^{\N}$ sets form an ideal. In fact we do not
know if $R^{\N} \cap {\cal H}$ is an ideal. We only have the following
result:
\begin{theorem}
  If $R^{\N} \cap {\cal H} $ is an ideal then $R^{\N} \cap {\cal H}$ is
  a $\sigma$-ideal.
\end{theorem}
\Proof
Suppose that $\{X_n : n \in \omega\}$ is an increasing sequence of
elements of $R^{\N} \cap {\cal H}$.
We will show that $X = \bigcup_{n \in \omega} X_n \in R^{\N}$.
We will use \ref{oldstuff}. Let $f \in \omega^\omega$ and suppose
that $x \leadsto S^x \in \Sigma_f$ is a Borel mapping.
By the assumption, for every $n$ there exists a function $g_n \in
{\cal X}_f$ such that 
$$\forall x \in X_n \ \forall^\infty k \ g_n(k) \not\in S^x(k).$$
As before, for $x \in X$ define $h^x \in \omega^\omega$ as follows:
$$h^x(n)=\min\{ k : \forall j > k \ g_n(j) \not \in S^x(j)\} \text{ for
  } n \in \omega.$$

Since mapping $x \leadsto h^x$ is Borel  we can find an increasing
function $h \in \omega^\omega$ such that $h^x \leq^\star  h$ for all $x \in
X$.
Let 
$$g(k) = g_n(k) \text{ for } h(n)<k<h(n+1).$$
It is easy to see that 
$$\forall x \in X \ \forall^\infty n \ g(n) \not \in S^x(n)$$
which finishes the proof.~$\QED$

\section{Strongly meager sets}
In this section we will estimate the smallest size of a set which is
not strongly meager. 

We start with the following forcing characterization of strongly
meager sets.

\begin{theorem}
  The following are equivalent:
  \begin{enumerate}
  \item $X$ is strongly meager,
  \item for every countable model $M \thinks \ZFCs$ there exists $s$
    such that $X+s \subseteq \R(M)$,
  \item for every measure zero set $H \subseteq 2^\omega \times
  2^\omega$ there exists $s$ such that 
$$\bigcup_{x \in X} (H)_{x+s} \neq 2^\omega .$$
  \end{enumerate}
\end{theorem}
\Proof
$(1) \rightarrow (3)$ Suppose that $H \subseteq 2^\omega \times
2^\omega$ has measure zero. Let $y$ be such that $(H)^y$ has measure
zero. Since $X$ is strongly meager there exists $s$ such that $(X+s)
\cap (H)^y = \emptyset$. In particular, $y \not \in \bigcup_{x \in X}
(H)_{x+s}$.

\vspace{0.1in}

$(3) \rightarrow (2)$ Let $G=2^\omega \setminus \R(M)$ and let $H$ be
such that $(H)_x = G+x$.
Find $$s \in 2^\omega \setminus \bigcup_{x \in X} (H)_{x}= 2^\omega
\setminus (X+G).$$
Clearly $s$ is the real we are looking for.

\vspace{0.1in}

$(2) \rightarrow (1)$
Suppose that $G$ is a null set. Let $M$ be a countable model such that
$G \in M$. Note that if $x+s \in \R(M)$ then $x+s \not\in G$.~$\QED$

\begin{definition}
  For $f \in \omega^\omega$ define the following cardinal invariants:

$$\cov(\Sigma_f) = \min\{|{\cal A}| : {\cal A} \subseteq \Sigma_f \ \&
\ \forall g \in {\cal X}_f \ \exists S \in {\cal A} \ \exists^\infty n \ g(n) \in S(n)\}$$
and 
$$\cof(\Sigma_f) = \min\{|{\cal A}| : {\cal A} \subseteq \Sigma_f \ \&
\ \forall g \in {\cal X}_f \ \exists S \in {\cal A} \ \forall^\infty n
\ g(n) \in S(n)\}.$$ 
Let $\unif(\cal {SM})$ be the least size of a set which is not strongly meager.
\end{definition}
We start with the following fact:
\begin{theorem}[Shelah, Brendle-Just]\label{shel1}
  There exists a measure zero set $H \subseteq 2^\omega$ such that for
  every set $Y$, if $Y+H$ has measure zero then there is $z$
  such that $Y+z \subseteq H$.
\end{theorem}
\Proof We shall describe how to modify the set constructed on  page
\pageref{shelahset} to get the required set $H$.

\begin{lemma}\label{trivlem1}
  There exists an increasing  function $f \in \omega^\omega$ and a
  sequence $\{U_n : n \in \omega\}$ such that for all $n \in \omega$:
  \begin{enumerate}
  \item $f(n+1) \geq 2^{f(0)+ \cdots + f(n) + n}$,
  \item $U_n \subseteq 2^{[f(n),f(n+1))}$,
  \item $$\mu\lft1([U_n]\rgt1) \geq
    1-\left(1-\frac{1}{2^n}\right)^{f(n)},$$
  \item for every $X \subseteq 2^{[f(n),f(n+1))}$, $|X| \leq n2^{n+1}$,
      $X+U_n \neq 2^{[f(n),f(n+1))}$.
  \end{enumerate}
\end{lemma}
\Proof
Suppose that $f(n)$ and $U_{n-1}$ have been chosen.
Let $k$ be such that 
$$\left(1-\frac{1}{2^k}\right)^{n2^{n+1}} \geq
1-\left(1-\frac{1}{2^n}\right)^{f(n)}$$ 
and
$$\left(1-\frac{1}{2^{n+1}}\right)\left(1-\frac{1}{2^k}\right)^{n2^{n+1}-1}\geq
\left(1-\frac{1}{2^{n}}\right).$$
Define $f(n+1)$ such that $f(n+1) \geq f(n)+kn2^{n+1} +n2^{f(n)}$.
Divide $[f(n),f(n+1))$ into $n2^{n+1}$ pairwise disjoint intervals of
size $k$, say $\{J_i : i \leq n2^{n+1}\}$.

Let 
$$U_n = \left\{s \in 2^{[f(n),f(n+1))}: \forall i \leq n2^{n+1} \
  \exists j \in J_i \ s(j)=1\right\} \text{ for } n \in \omega.$$
By the choice of $k$ we have
$$\mu\lft1([U_n]\rgt1) \geq
    1-\left(1-\frac{1}{2^n}\right)^{f(n)}.$$
Suppose that $X \subseteq 2^{[f(n),f(n+1))}$, $|X| \leq n2^{n+1}$. Let
  $X=\{x_i : i\leq n2^{n+1}\}$.
Define $x^\star$ such that $x^\star \rest J_i=x_i \rest J_i$ for every $i$.

It is easy to see that $x^\star \not \in U_n+X$.~$\QED$

As in \ref{trivlem}, define sets
$\{I_s : s \in 2^{f(n)}\}$ such that 
\begin{enumerate}
\item $I_s \subseteq \lft1[f(n), f(n+1)\rgt1)$ for all $s$,
\item $I_s \cap I_t = \emptyset$ for $s \neq t$,
\item $|I_s|=n-1$ for all $s$,
\item for every $s \in 2^{f(n)}$ there exists $i \leq n2^{n+1}$ such
  that $I_s \subseteq J_i$.
\end{enumerate}
For $s \in 2^{f(n)}$ let 
$$A_s =\left\{t \in U_n : \exists k \in I_s \
t(k)=1\right\}.$$ 

Note  that   the choice of $k$ and (1)-(4) above guarantee that  the sets
$A_s$ satisfy the conditions of 
\ref{trivlem} (except that the measure of $[A_s]$ is different). 

Let $T_m$'s and $H$ be defined as in \ref{charnuladd}. 

Let $\widehat{T}=\prod_{n \in \omega} U_n$. Note that $\bigcup_m [T_m]
\subseteq \rationals + \widehat{T}$.

As in the proof of \ref{charnuladd}, we show that if $Y$ is a 
set such that $Y+H$ has measure zero then
there exists 
  a sequence $\{I_n : n \in \omega\}$ such that 
  \begin{enumerate}
  \item for all $n$, $I_n \subseteq 2^{[f(n), f(n+1))}$,
  \item for all $n$, $|I_n| \leq n2^{n+1}$,
  \item $\forall x \in Y \ \forall^\infty n \ x \rest \lft1[f(n),
    f(n+1)\rgt1) \in I_n$.
\end{enumerate}

Note that if $Y$ has the above property then $Y + \widehat{T}+\rationals
\neq 2^\omega$. This is guaranteed by \ref{trivlem1}(4). 
Let $z \in 2^\omega \setminus (Y + \widehat{T}+\rationals)$. 
Clearly 
$$z+ Y \subseteq 2^\omega \setminus (\widehat{T}+\rationals)
\subseteq 2^\omega \setminus \bigcup_m [T_m] =H,$$
which finishes the proof.~$\QED$

Now we can characterize strongly meager sets.

\begin{theorem}
  $\min_f \cov(\Sigma_f) \leq \unif({\cal {SM}}) \leq \min_f \cof(\Sigma_f).$
\end{theorem}
\Proof The first inequality is proved like \ref{oldstuff} (see
\cite{Bar88Cov} for details).

\vspace{0.1in}

Fix $f \in \omega^\omega$. Let ${\cal A} \subseteq
\Sigma_f$ be a family of size $\cof(\Sigma_f)$ such that 
$$\forall g \in {\cal X}_f \ \exists S \in {\cal A} \ \forall^\infty n
\ g(n) \in S(n).$$ 
By increasing $f$ and doing some elementary coding 
we can assume that there is a function
$\widetilde{f}$ such that $f(n) = 2^{\tilde{f}(n+1)-\tilde{f}(n)}$ for
all $n$. 
For $S \in \Sigma_f$ define 
$$P_S = \left\{x \in 2^\omega: \forall^\infty n \ x \rest
\lft1[\widetilde{f}(n), \widetilde{f}(n+1)\rgt1) \in S(n)\right\},$$
where $S(n)$ is treated as a subset of $2^{[\tilde{f}(n),
  \tilde{f}(n+1))}$.
Let $H $ be the set from \ref{shel1}. 

Note that \ref{trivlem1} gives us enough freedom to ensure 
that $P_S + H$ has measure zero for all $S \in \Sigma_f$.
This is because in \ref{trivlem1}(4) we can replace the clause $|X|<n
2^{n+1}$ by 
$|X|<g(n)$, where $g$ is an arbitrary fixed function. (cf.
\cite{ErdKunMaul81Som}, lemma 9).

For each $S \in {\cal A}$ choose $x_S \in 2^\omega$ such that $P_S +
x_S \subseteq H$. 
We claim that $X=\{x_S : S \in {\cal A}\}$ is not strongly meager and 
 show that  $X+H = 2^\omega$.
Suppose that $y \in 2^\omega \setminus (X+H)$. Let $S \in {\cal A}$ be such
that $y \in P_S$. Now $y \in P_S \subseteq H+x_S \subseteq H+X$.
Contradiction.~$\QED$ 

A forcing notion ${\cal P}$ has Laver property if it does not
``increase'' $\cof(\Sigma_f)$. In other words, if ${\cal P}$ has Laver
property then
for every $g \in
\V^{{\cal P}} \cap {\cal X}_f$ there exists $S \in \V \cap \Sigma_f$
such that $g(n) \in S(n)$ for all $n$.

In this language we can formulate the previous theorem as follows:

\begin{theorem}
  Suppose that $\V'$ is a generic extension of $\V$ obtained by a
  forcing notion which has Laver property. Then the set $\V \cap
  2^\omega$ is not strongly meager in $\V'$.~$\QED$
\end{theorem}

\vspace{0.1in}

{\bf Acknowledgments}

We would like to thank Joerg Brendle for reading the manuscript and
for his helpful remarks. We would also like to thank the referees of
this paper for shaping  the final version of this paper.

\ifx\undefined\bysame
\newcommand{\bysame}{\leavevmode\hbox to3em{\hrulefill}\,}
\fi


\begin{thebibliography}{10}

\bibitem{Bar84Add}
Tomek Bartoszy\'{n}ski, {\em Additivity of measure implies additivity of
  category}, Transactions of the American Mathematical Society {\bf 281}
  (1984), no.~1, 209--213.

\bibitem{Bar88Cov}
\bysame, {\em On covering of real line by null sets}, Pacific Journal of
  Mathematics {\bf 131} (1988), no.~1, 1--12.

\bibitem{Bar87Com}
\bysame, {\em Combinatorial aspects of measure and category}, Fundamenta
  Mathematicae {\bf 127} (1987), no.~3, 225--239.

\bibitem{BarJud93Cov}
Tomek Bartoszy\'{n}ski and Haim Judah, {\em On the smallest covering of the
  real line by meager sets {II}}, Proceedings of the American Mathematical
  Society, to appear.

\bibitem{BJS}
Tomek Bartoszy\'{n}ski, Winfried Just, and Marion Scheepers, {\em Covering
  games and the {B}anach-{M}azur game: k-tactics}, The Canadian Journal of
  Mathematics {\bf 45} (1993), no.~5, 897--929.

\bibitem{ErdKunMaul81Som}
Paul Erdos, K.~Kunen, and R.~Daniel Mauldin, {\em Some additive properties of
  sets of real numbers}, Fundamenta Mathematicae {\bf
  113} (1981), no.~3, 187--199.

\bibitem{FrePa}
David~H. Fremlin, {\em The partial orderings in measure theory and {T}ukey
  ordering}, Note di Matematica, to appear.

\bibitem{HuSeq27}
Witold Hurewicz, {\em On sequences of continuous functions}, Fundamenta
  Mathematicae {\bf 9} (1927), 193--204.

\bibitem{Kur66Top}
K.~Kuratowski, {\em Topology}, New edition, revised and augmented. Translated
  from the French by J. Jaworowski, vol.~I, London; Panstwowe Wydawnictwo
  Naukowe, Warsaw, 1966.

\bibitem{Mil81Som}
Arnold~W. Miller, {\em Some properties of measure and category}, Transactions
  of the American Mathematical Society {\bf 266} (1981), no.~1, 93--114.

\bibitem{Mil84Add}
\bysame, {\em Additivity of measure implies dominating reals}, Proceedings of
  the American Mathematical Society {\bf 91} (1984), no.~1, 111--117.

\bibitem{FreMil88Som}
Arnold~W. Miller and David~H. Fremlin, {\em On some properties of {H}urewicz,
  {M}enger, and {R}othberger}, Fundamenta Mathematicae
  {\bf 129} (1988), no.~1, 17--33.

\bibitem{Paw85Pow}
Janusz Pawlikowski, {\em Powers of transitive bases of measure and category},
  Proceedings of the American Mathematical Society {\bf 93} (1985), no.~4,
  719--729.

\bibitem{Paw92Sie}
\bysame, {\em All {S}ierpi\'{n}ski sets are strongly meager}, 1992.

\bibitem{PawPropC}
\bysame, {\em Property ${C}''$, strongly meager sets and subsets of the plane},
  1993.

\bibitem{PawRecPar}
Janusz Pawlikowski and Ireneusz Rec{\l}aw, {\em Parametrized {C}ichon's diagram
  and small sets}, 1993.

\bibitem{RecOpen}
Ireneusz Rec{\l}aw, {\em Every {L}usin set is undetermined in point--open
  game}, Fundamenta Mathematicae, {\bf 144}, (1994), 43--54


\bibitem{Rec89Sma}
\bysame, {\em On small sets in the sense of measure and category}, Polska
  Akademia Nauk. Fundamenta Mathematicae {\bf 133} (1989), no.~3, 255--260.

\bibitem{Sh445}
Saharon Shelah, {\em Every null additive set is meagre additive}, Israel
  Journal of Mathematics, {\bf 89}, 357--376, publication 445.

\bibitem{vdo:inttop}
E.~K. van Douwen, {\em The integers and topology}, Handbook of Set Theoretic
  Topology (Amsterdam) (K.~Kunen and J.~E. Vaughan, eds.), North-Holland,
  Amsterdam, 1984, pp.~111--167.
\end{thebibliography}
\end{document}